\documentclass[12pt,leqno]{article}
\usepackage{amssymb}
\usepackage{amsmath}
\usepackage{enumerate}
\usepackage[latin1]{inputenc}
\def\S{\mathbb{S}}
\def\H{\mathbb{H}}
\def\R{\mathbb{R}}
\def\C{\mathbb{C}}
\def\P{\mathbb{P}}

\newtheorem{remark}{Remark} 
\newtheorem{theorem}{Theorem}
\newtheorem{corollary}{Corollary}

\begin{document}

\title{Lagrangian surfaces with circular ellipse of curvature in complex space forms}         
\author{Ildefonso Castro\thanks{Research partially supported by a MCYT grant No. BFM2001-2967.}}        


\maketitle

\pagestyle{plain}

\begin{abstract}
We classify the Lagrangian orientable surfaces in complex space forms with the property that the ellipse of curvature is always a circle. As a consequence, we obtain new characterizations of the Clifford torus in the complex projective plane and of the Whitney spheres in the complex projective, complex Euclidean and complex hyperbolic planes.


\end{abstract}

\section{Introduction}       

An interesting notion that comes up in the study of surfaces in higher codimension is that of the ellipse of curvature. This is the image in the normal space of the unit circle in the tangent plane under the second fundamental form.

Using this concept, I.V. Guadalupe and L. Rodriguez obtained in [GR] some inequalities relating the area of compact surfaces in (real) space forms and the integral of the square of the norm of the mean curvature vector with topological invariants. When the ellipse of curvature is a circle, restrictions on the Gaussian and normal curvatures gave them some rigidity results. (See Remark 2 in the paper).

M.F. Atiyah and H.B. Lawson showed (see Remark 1 in [GR]) that an immersed surface in the 4-sphere $\S^4$ has the ellipse always a circle  if and only if the canonical lift of the immersion into the bundle of almost complex structures of $\S^4$ is holomorphic. Holomorphic curves in this twistor space of $\S^4$ can also be projected down to $\S^4$ in order to obtain examples of surfaces in $\S^4$ verifying that the ellipse of curvature is always a circle.

Other works which use the ellipse of curvature as a tool in the study of surfaces in (real) space forms can be found in [MRR], [R], [V] and [W].

\vspace{0.1cm}

In this article we attach the circular ellipse of curvature condition to Lagrangian surfaces in complex space forms. Here a complex space form means the complex projective plane $\C\P^2$, the complex Euclidean plane $\C^2$ or the complex hyperbolic plane $\C\H^2$, endowed with their canonical structures of Kaehler surfaces. An isometric immersion of a surface in a complex space form $N(c)$ is
Lagrangian if the complex structure of $N(c)$ carries each tangent plane into the corresponding normal plane.

There exist many interesting Lagrangian surfaces in complex space forms.
F. Urbano and the author recovered a classical example in Symplectic Geometry, known as the {\em Whitney sphere}, given (up to translations and scaling) by
\begin{equation}
\label{eq:WhitneyC2}
\Phi: \S^2 \longrightarrow \C^2, \quad \Phi (x,y,z) = \frac{1+i z}{1+z^{2}} 
   \left(  x, y \right),
\end{equation}
showing that, from a Riemannian point of view, this Lagrangian sphere plays the role of the round sphere in the Lagrangian setting (see [CU1]).

Also, in [CU3], F. Urbano and the author completely classified all the twistor holomorphic Lagrangian immersions in the complex projective plane $\C\P^2$
(resp. in the complex hyperbolic plane $\C\H^2$), i.e. those Lagrangian immersions such that their twistor lifts to the twistor space over $\C\P^2$
(resp. over $\C\H^2$) are holomorphic. (Later the classification of the twistor harmonic ones was obtained in [CU4]).

The first classification provided a 
one-parameter family of examples of Lagrangian spheres in $\C\P^2$, 
given by
\begin{equation}
\label{eq:WhitneyCP2}
\begin{array}{c}
 \Phi_t : \S^2 \longrightarrow \C\P^2 , \, \, t\geq  0, \\
  \Phi_t(x,y,z) =\left[
    \frac{1}{c_t^2+s_t^2 z^2}   
    \left( c_t x, s_t xz, c_t y, s_t yz, z,
    s_t c_t (1+z^2) \right) \right], 
\end{array}
\end{equation}
(with $c_t=\cosh t $, $s_t = \sinh t$). We observe that $\Phi_0$ is the totally geodesic Lagrangian immersion of $\S^2$ in $\C\P^2$.
In $\C\H^2$, we got a bigger family of examples:
\begin{equation}
\label{eq:WhitneyCH2}
\begin{array}{c}
\widehat{\Phi}_t : \S^2 \longrightarrow \C\H^2 , \,\, t >0, \\
  \widehat{\Phi}_t(x,y,z) =\left[ \frac{1}{s_t^2+c_t^2z^2}
        \left( s_t x, c_t xz, s_t y, c_t yz, z,
    -s_t c_t (1+z^2) \right)\right]; 
\end{array}
\end{equation}

\begin{equation}
\label{eq:PhisCH2}
\begin{array}{c}
 \Psi_s : \C ^* \longrightarrow \C\H^2 , \, \, s\in [0,\pi /4[, \\
 \Psi_s(z)=\left[ \frac{1}{|w(z)|^2}
   \left( \cos ^2 \! s \,z^2 - \sin ^2 \! s \,\overline{z}^2,
     \frac{|z|^2-1}{\sqrt 2} w(z),
     \frac{|z|^2+1}{\sqrt 2} w(z)
   \right) \right], 
\end{array}
\end{equation}
where $w(z)=\cos s \,z + \sin s \,\overline{z}$; and
\begin{equation}
\label{eq:etaCH2}
\begin{array}{c}
  \eta:\C \longrightarrow \C\H^2 , \, \, \, \eta(x,y)=  \\
  \left[ \frac{1}{1+4x^2}\left( 2y,4xy,2x-4x^3-4xy^2,6x^2+2y^2, 1+6x^2+2y^2,
                          4x^3+4xy^2 \right) \right]. 
\end{array}
\end{equation}

It is a well known fact that Lagrangian immersions in $\C\P^2$ (resp. in $\C\H^2$) are ---locally--- projections by the Hopf fibration $\Pi: \S^5 \longrightarrow \C\P^2$  (resp. $\Pi: \H^5_1 \longrightarrow \C\H^2$) of horizontal immersions in the unit sphere $\S^5$ (resp. in the anti-De Sitter space $\H^5_1$).
Above $[\,\, ]$ means image under these Hopf fibrations.

Both families (\ref{eq:WhitneyCP2}) and (\ref{eq:WhitneyCH2}) are also
called {\em Whitney spheres} in [CMU] because of their geometric properties, similar to the classical Whitney spheres in $\C^2$. Examples (\ref{eq:PhisCH2}) and (\ref{eq:etaCH2}) are complete Lagrangian embeddings with zero total curvature. All these examples 
(\ref{eq:WhitneyCP2})--(\ref{eq:etaCH2}) were studied by B.Y. Chen and L. Vrancken in a different context (see [CV]).

\vspace{0.3cm}

Making use of the mentioned results and after detecting a simple close relationship between the ellipse of curvature and some differential forms naturally associated to a Lagrangian surface, we prove 
in Theorem 1 that 
{\em the only Lagrangian orientable surfaces with circular ellipse of curvature in complex space forms are the minimal ones and (open subsets of) the examples (\ref{eq:WhitneyC2}), (\ref{eq:WhitneyCP2}), (\ref{eq:WhitneyCH2}), (\ref{eq:PhisCH2}) and (\ref{eq:etaCH2})}.

Among the main consequences of this Theorem we emphasize the following new characterizations (Corollary 1) of the examples (\ref{eq:WhitneyC2}), (\ref{eq:WhitneyCP2}) and (\ref{eq:WhitneyCH2}):
\begin{quote}
{\em 
\begin{enumerate}
\item The Whitney spheres $\Phi: \S^2 \longrightarrow \C^2$ are the only Lagrangian compact orientable surfaces with circular ellipse of curvature in complex Euclidean plane.
\item The Whitney spheres $\Phi_t : \S^2 \longrightarrow \C\P^2 , \, \, t> 0,$ are the only non minimal
Lagrangian compact orientable surfaces with circular ellipse of curvature in the complex projective plane.
\item The Whitney spheres $\widehat{\Phi}_t : \S^2 \longrightarrow \C\H^2 , \,\, t >0,$ are the only Lagrangian compact orientable surfaces with circular ellipse of curvature in the complex hyperbolic plane.
\end{enumerate}
}
\end{quote}
In $\C\P ^2$, we also determine the Clifford torus in terms of the circular ellipse of curvature condition and a pinching on its radius function (Corollary 2).

\section{Preliminaries}

Let $\phi : \Sigma \longrightarrow N $ be an immersion of a surface $\Sigma $ into a 4-dimensional Riemannian manifold $N$. We denote the metric of $N$ as well as the induced metric in $\Sigma $ by $\langle, \rangle$. If $\sigma: T\Sigma \times T\Sigma \longrightarrow T^\perp \Sigma $ is the second fundamental form of $\phi$, the {\em ellipse of curvature} is the subset of the normal plane defined as $\{ \sigma (v,v) \in T_p^\perp \Sigma \,:\, \langle v,v \rangle =1,
v\in T_p\Sigma, \, p\in \Sigma \}$. To see that it is an ellipse, we consider an arbitrary orthogonal tangent frame $\{ e_1, e_ 2 \}$, denote $\sigma_{ij}=\sigma (e_i,e_j)$, $i,j=1,2$, and look at the following formula for $v=\cos \theta \, e_1 + \sin \theta \, e_2$:
\begin{equation}
\label{eq:ellipse}
\sigma (v,v) = H + \cos 2 \theta \, \frac{\sigma_{11}-\sigma_{22}}{2} 
+ \sin 2 \theta \, \sigma_{12}, 
\end{equation}
where $H= 1/2$ trace $\sigma$ is the mean curvature vector. 
As $v$ goes once around the unit tangent circle, $\sigma (v,v)$ goes twice around the ellipse, which could degenerate into a line segment or a point.
The center of the ellipse is $H$. $\phi $ is said to be minimal if $H\equiv 0$.

From (\ref{eq:ellipse}) it is not difficult to deduce that the ellipse of curvature is a circle if and only if
\begin{equation}
\label{eq:circle}
\frac{|\sigma_{11}-\sigma_{22}|^2}{4} = |\sigma_{12}|^2, \, \, 
\langle \sigma_{11}-\sigma_{22}, \sigma_{12} \rangle =0.
\end{equation}
The property that the ellipse is a circle is a conformal invariant.
\vspace{0.5cm}

We consider now that the target manifold $N$ is a simply connected complex space form with complex structure $J$ and constant holomorphic sectional curvature $c$. We denote by $N(c)$ the complex projective plane $\C\P^2$ if $c=4$, the complex Euclidean plane $\C^2$ if $c=0$ and the complex hyperbolic plane $\C\H^2$ if $c=-4$, with their standard complex structures and metrics. An immersion $\phi: \Sigma \longrightarrow N(c) $ is said to be {\em Lagrangian} (it is also used classically the nomenclature of {\em totally real} instead of Lagrangian) if $JT_p\Sigma = T_p^\perp \Sigma$
for each point $p$ of $M$. In this case it is not hard to find that
the trilinear form 
\[
C(u,v,w)=\langle \sigma (u,v), J\phi_*(w) \rangle 
\]
is fully symmetric and, according to the Gauss equation of $\phi$, 
$K=c/4+\langle \sigma_{11}, \sigma_{22} \rangle - |\sigma_{12}|^2$,
the Gauss curvature $K$ of a Lagrangian surface in $N(c)$ is given by
\begin{equation}
\label{eq:Gauss}
K=\frac{c}{4}+2|H|^2-\frac{|\sigma|^2}{2}.
\end{equation}
If the ellipse of curvature of a surface is always a circle, we can define its radius function $R\geq 0$ by means of
$R^2=(|\sigma_{11}-\sigma_{22}|^2)/4 = |\sigma_{12}|^2$ 
taking into account (\ref{eq:circle}). From (\ref{eq:Gauss}), it is not complicated to obtain that the radius function of a Lagrangian surface with circular ellipse of curvature in a complex space form $N(c)$ is given by
\begin{equation}
\label{eq:R}
R^2= \frac{1}{2}\left(\frac{c}{4}+|H|^2-K \right).
\end{equation}

\section{Statements and proofs of the results}

\begin{theorem}
Let $\phi : \Sigma \longrightarrow N(c)$ be a Lagrangian immersion of an orientable surface $\Sigma $ into a complex space form $N(c)$. The ellipse of curvature of $\phi $ is always a circle if and only if
either $\phi $ is minimal or $\phi $ is locally congruent to some of the immersions (\ref{eq:WhitneyC2}), (\ref{eq:WhitneyCP2}), (\ref{eq:WhitneyCH2}), (\ref{eq:PhisCH2}) or (\ref{eq:etaCH2}). 
\end{theorem}

\begin{remark}
{\rm The ellipse of curvature in a totally geodesic point degenerates into a point. This case can be thought as a circular ellipse of curvature with radius equal to zero.}
\end{remark}

{\it Proof:\/} 
We start defining  a cubic differential form $\Theta $ and a differential form $\Upsilon $ (see [CU1] and [CU3]) by 
\[
\begin{array}{c}
\Theta (z) = f(z)(dz)^3, \, f=4\, C(\partial_z,\partial_z,\partial_z), \\ \\
\Upsilon (z) = h(z)\,dz, \, h=2 \, C(\partial_z,\partial_{\overline z},\partial_{\overline z}) / \langle \partial_z, \partial_{\overline z}\rangle ,
\end{array}
\]
where $z=x+iy$ is a local isothermal coordinate on $\Sigma $ where the induced metric is written as $\langle , \rangle = e^{2u} |dz|^2$ and
$\langle,\rangle$, $J$ and $\sigma $ are extended $\C$-linearly to the complexified bundles. We note that $2\sigma(\partial_z,\partial_{\overline z})=e^{2u}H$.

We choose the orthonormal frame 
$\{ e_1=\partial_x/e^u,e_2=\partial_y/e^u \} $ and put $C_{ijk}=C(e_i,e_j,e_k)=\langle \sigma(e_i,e_j),J\phi_* e_k \rangle$,
$i,j,k\in \{ 1,2 \}$.
We can now rewrite $f$ and $h$ as follows:
\begin{equation}
\label{eq:fh}
\begin{array}{c}
f=\frac{e^{3u}}{2}[(C_{111}-3\,C_{122})+i\,(C_{222}-3\,C_{112})], \\ \\
h=\frac{e^{u}}{2}[(C_{111}+C_{122})+i\,(C_{112}+C_{222})].
\end{array}
\end{equation}
On the other hand, using again the symmetry of $C$, we obtain:
\begin{equation}
\label{eq:q}
\begin{array}{c}
\left| \sigma_{11}-\sigma_{22} \right|^2- 4\,|\sigma_{12}|^2= 
 C_{111}^2 - 2\, C_{111}\,C_{112}-\\ \\ 
-3\,C_{122}^2 - 3 \,C_{112}^2 -2 \,C_{112} \, C_{222}+C_{222}^2 ,\\ \\
\langle \sigma_{11}-\sigma_{22}, \sigma_{12} \rangle = 
C_{111}\, C_{112}-C_{122}\, C_{222}.
\end{array}
\end{equation}
From (\ref{eq:fh}) and (\ref{eq:q}) we get:
\begin{equation}
\label{eq:formula}
e^{-4u}f\, \overline{h} = 
\left( \frac{|\sigma_{11}-\sigma_{22}|^2}{4} - |\sigma_{12}|^2 \right)
+i  \,
\langle \sigma_{11}-\sigma_{22}, \sigma_{12} \rangle .
\end{equation}
Using (\ref{eq:circle}), we know that ellipse of curvature is a circle if and only if the right-side of the equality in (\ref{eq:formula}) vanishes.

In [CU1, Theorem 2] and [CU3, Theorems 1 and 2] it was proved that the cubic differential $\Theta $ is identically zero (i.e. $f\equiv 0$) if and only if $\phi $ is totally geodesic or $\phi (\Sigma)$ is an open set of $\Phi(\S^2)$ (see (\ref{eq:WhitneyC2})) if $c=0$, of $\Phi_t(\S^2)$ (see (\ref{eq:WhitneyCP2})) if $c=4$, of $\Phi_t(\S^2)$ (see (\ref{eq:WhitneyCH2})), $\Psi_s(\C^*)$ (see (\ref{eq:PhisCH2})) or $\eta (\C)$
(see (\ref{eq:etaCH2})) if $c=-4$. Moreover, $\phi $ is minimal if and only if $h\equiv 0$ since $|h|^2 = e^{2u} |H|^2$. 

A standard argument completes the proof of the theorem.$_\diamondsuit$
\vspace{0.1cm}

Minimal Lagrangian orientable surfaces in $\C^2$ can be represented as holomorphic curves in $\C^2$ by exchanging the orthogonal complex structure on $\R^4$ (cf. [CM]). In the case of $\C\P^2$ we emphasize the following result.
\begin{theorem}[\mbox{[Y], Theorem 7; [NT]}]
Let $\phi :\Sigma \longrightarrow  \C\P ^2$ be a minimal 
Lagrangian immersion of a compact orientable surface. Then:
\begin{description}
\item[{\rm a)}] If  $\Sigma $ has genus zero, then $\phi $ is totally geodesic
and it is the standard immersion $\Phi_0$ (see (\ref{eq:WhitneyCP2})) of $\S^2$ into $\C \P ^2$.
\item[{\rm b)}] If the genus of $\Sigma $ is nonzero and the Gauss curvature $K$ is either nonnegative or nonpositive, then $\Sigma $ is flat and $\phi $ is the Clifford torus.
\end{description}
\end{theorem}

The Clifford torus is the quotient by the Hopf fibration 
of the isometric inclusion
\[ \S^{1}(1/\sqrt3) \times \S^{1}(1/\sqrt3) \times \S^{1}(1/\sqrt3)
 \hookrightarrow \S^{5} \subset \C^3. \] 
It was the only example of minimal Lagrangian torus 
in $\C\P^2 $ during a long period. For subsequent progress in this matter we refer to [CU2] and [MM].

Using the well known fact that there do not exist compact minimal surfaces neither in $\C^2$ nor in $\C\H^2$, we deduce from Theorem 1 the following first immediate consequence.

\begin{corollary}
The Whitney spheres (\ref{eq:WhitneyC2}) (resp. (\ref{eq:WhitneyCH2})) are the only Lagrangian compact orientable surfaces with circular ellipse of curvature in the complex Euclidean plane (resp. in the complex hyperbolic plane).
The Whitney spheres (\ref{eq:WhitneyCP2}) are the only non minimal Lagrangian orientable compact  surfaces with circular ellipse of curvature in the complex projective plane.
\end{corollary}

\begin{remark}
{\rm  In [CU1] and [CU3] we obtained a different characterization of the Whitney spheres that can be written as follows: 

{\em Let $\phi : \Sigma \longrightarrow N(c)$ be a Lagrangian
immersion of a compact orientable surface $\Sigma $ with Euler characteristic $\chi (\Sigma)$.  Then:
\[
 \int_{\Sigma} |H|^2 dA + \frac{c}{2}  {\rm Area}\, (\Sigma)
\geq 4 \pi  \chi (\Sigma) ,
\]
and the equality holds if and only if $\phi $ is congruent to some of the
Whitney spheres.}

It is interesting to compare this result with Theorem 1 in [GR] since both
follows a similar style.}
\end{remark}

\vspace{0.2cm}

On the other hand, O. Kasabov proved in [K] the following result: {\em
Let $\phi : \Sigma \longrightarrow N(c)$ be a constant $\lambda$-isotropic Lagrangian immersion of a complete surface $\Sigma $ into a complex space form $N(c)$. Then either $\phi $ is totally geodesic or $c>0$ ($c=4$) and $\phi $ is congruent to the Clifford torus.} 

The Lagrangian immersion $\phi $ is said to be $\lambda$-isotropic ($\lambda$ being a constant) if for each unit vector $v$ tangent to $\Sigma $ the length of $\sigma (v,v)$ is $\lambda$. In other words, $\phi $ must be a minimal immersion whose radius function (of the circular ellipse of curvature, see (\ref{eq:R})) is $\lambda $=constant. For example, the Clifford torus is $1/\sqrt{2}$-isotropic in $\C\P^2$.
We for our part finish the paper showing the following alternative determination of the Clifford torus.

\begin{corollary}
Let $\phi : \Sigma \longrightarrow N(c)$ be a Lagrangian immersion with circular ellipse of curvature and radius function $R$ of a compact orientable surface $\Sigma $ into a complex space form $N(c)$. 

If $R\geq 1/\sqrt{2}$, then $c=4$ and $\phi $ is congruent to the Clifford torus.
\end{corollary}

{\it Proof:\/} From Theorem 1 and the compactness hypothesis, our immersion $\phi $ must be minimal (necessarily in $\C\P^2$) or congruent to some of the examples (\ref{eq:WhitneyC2}), (\ref{eq:WhitneyCP2}) or (\ref{eq:WhitneyCH2}). 

Using (\ref{eq:Gauss}) and (\ref{eq:R}), if $\phi $ is minimal its radius function $R$ is given by $R^2=1/2(1-K)=|\sigma|^2/4$. But $R\geq 1/\sqrt{2}$ implies $K\leq 0$ and Theorem 2 says to us that $\phi $ is the Clifford torus.

Otherwise, the cubic differential $\Theta $ associated to a Whitney sphere
vanishes identically (see the proof of Theorem 1). Since
$|f|^2=e^{6u}(c/2+|H|^2-2K)\equiv 0$, we deduce that $R^2=1/2(K-c/4)$ in this case. Now the hypothesis $R\geq 1/\sqrt{2}$ translates into $K\geq 1$
for (\ref{eq:WhitneyC2}), $K\geq 2$ for (\ref{eq:WhitneyCP2}) and $K\geq 0$
for (\ref{eq:WhitneyCH2}). But in [CU1] and [CU3] we studied the behavior of the Gauss curvature of the Whitney spheres showing:
\begin{description}
\item[{\rm (i)}]  The Gauss curvature $K$ of $(\S^2,\Phi ^* \langle, \rangle)$ verifies
$0\leq K \leq 1$. Moreover, $K(x,y,z)=0$ (respectively  $K(x,y,z)=1$)
if and only if $z=\pm 1$ (respectively $z=0$).
\item[{\rm (ii)}] The Gauss curvature $K_t$ of $(\S^2,
\Phi_t ^* \langle, \rangle)$ verifies
 $1 \leq K_t \leq 1+2\sinh ^2 t$. Moreover, $K_t(x,y,z)=1$
 (respectively  $K_t(x,y,z)=1+2\sinh ^2 t$)
 if and only if $z=\pm 1$ (respectively $z=0$).
\item[{\rm (iii)}] The Gauss curvature $K_t$ of 
   $(S^2,\widehat{\Phi}_t^* \langle, \rangle)$
   verifies  $-1 \leq K_t \leq -1+2\cosh ^2 t$. Moreover, $K_t(x,y,z)=-1$
   (respectively $K_t(x,y,z)=-1+2\cosh ^2 t$)
   if and only if $z=\pm 1$ (respectively $z=0$).
\end{description}
This is incompatible with our previous conclusion. So the second case is impossible.$_\diamondsuit$

\vspace{1cm}

\begin{tabular}{l}
{\sc address}:  \\
Departamento de Matem\'{a}ticas   \\
Escuela Polit\'{e}cnica Superior  \\
Universidad de Ja\'{e}n \\
23071 Ja\'{e}n \\
SPAIN \\
{\sc e-mail}:  \\
{\tt icastro@ujaen.es} \\
\end{tabular}

\end{document}